\documentclass[a4paper,12pt,reqno]{amsart}
\usepackage{lmodern}
\usepackage{graphicx,amsmath,amssymb,amsthm}
\numberwithin{equation}{section}
\newtheorem{Theorem}{Theorem}[section]

\newtheorem{Corollary}{Corollary}[section]

\theoremstyle{remark}
\newtheorem{Remark}{Remark}[section]
\newtheorem{Example}{Example}[section]
\title[]{On the behavior of solutions of quasilinear elliptic inequalities near a boundary point}
\author{Andrej A. Kon'kov}
\address{Department of Differential Equations,
Faculty of Mechanics and Mathematics,
Mo\-s\-cow Lo\-mo\-no\-sov State University,
Vorobyovy Gory,
Moscow, 119992 Russia}
\email{konkov@mech.math.msu.su}
\date{}

\begin{document}

\begin{abstract}
Assume that $p > 1$ and $p - 1 \le \alpha \le p$ are real numbers and
$\Omega$ is a non-empty open subset of ${\mathbb R}^n$, $n \ge 2$.
We consider the inequality
$$
	{\rm div} \, A (x, D u)
	+
	b (x) |D u|^\alpha
	\ge
	0,
$$
where 
$D = (\partial / \partial x_1, \ldots, \partial / \partial x_n)$
is the gradient operator and
$A : \Omega \times {\mathbb R}^n \to {\mathbb R}^n$ 
and 
$b : \Omega \to [0, \infty)$ 
are some functions with
$$
	C_1
	|\xi|^p
	\le
	\xi
	A (x, \xi),
	\quad
	|A (x, \xi)|
	\le
	C_2
	|\xi|^{p-1}, 
	\quad
	C_1, C_2 = const > 0,
$$
for almost all
$x \in \Omega$
and for all
$\xi \in {\mathbb R}^n$.
For solutions of this inequality we obtain estimates depending on the geometry of $\Omega$.
In particular, these estimates imply regularity conditions of a boundary point.
\end{abstract}

\maketitle

\section{Introduction}
Let $\Omega$ be an open subset of ${\mathbb R}^n$, $n \ge 2$.
By $B_r^x$ and $S_r^x$ we mean the open ball and the sphere in ${\mathbb R}^n$ 
of radius $r > 0$ and center at a point $x$.
In the case of $x = 0$, we write $B_r$ and $S_r$
instead of $B_r^0$ and $S_r^0$, respectively.
Let us denote
$
	B_{r_1, r_2}
	=
	\{
		x \in {\mathbb R}^n : r_1 < |x| < r_2
	\}
$
and
$
	\Omega_{r_1, r_2}
	=
	B_{r_1, r_2}
	\cap
	\Omega,
$
$0 < r_1 < r_2$.
Through out the paper, we assume that $S_r \cap \Omega \ne \emptyset$ for any $r \in (0, R)$,
where $R > 0$ is some real number.

We are interested in the behavior of solutions of the problem
\begin{equation}
	{\rm div} \, A (x, D u)
	+
	b (x) |D u|^\alpha
	\ge
	0
	\quad
	\mbox{in } B_R \cap \Omega,
	\quad
	\left.
		u
	\right|_{
		B_R \cap \partial \Omega
	}
	=
	0,
	\label{1.1}
\end{equation}
where 
$D = (\partial / \partial x_1, \partial / \partial x_2, \ldots, \partial / \partial x_n)$
is the gradient operator and the function
$A : \Omega \times {\mathbb R}^n \to {\mathbb R}^n$ 
satisfies the ellipticity conditions
$$
	C_1
	|\xi|^p
	\le
	\xi
	A (x, \xi),
	\quad
	|A (x, \xi)|
	\le
	C_2
	|\xi|^{p-1}
$$
with some constants $C_1 > 0$, $C_2 > 0$, and $p > 1$
for almost all $x \in \Omega$ and for all $\xi \in {\mathbb R}^n$.
It is also assumed that $p - 1 \le \alpha \le p$ is a real number and 
$b$ is a non-negative function such that $b \in L_\nu (\Omega_{r, R})$ for all $r \in (0, R)$,
where $\nu$ satisfies the following requirements:
\begin{enumerate}
\item[(i)]
if $\alpha = p$, then $\nu = \infty$;
\item[(ii)] 
if $\alpha = p - 1$ and $n \ne p$, then $\nu = \max \{ n, p \}$;
\item[(iii)] 
if $\alpha = p - 1$ and $n = p$, then $\nu > p$;
\item[(iv)] 
if $p - 1 < \alpha < p$ and $n \ne p$, then $\nu = \max \{ n, p \} / (p - \alpha)$;
\item[(v)] 
if $p - 1 < \alpha < p$ and $n = p$, then $\nu > p / (p - \alpha)$.
\end{enumerate}

We say that
$
	u 
	\in 
	{W_p^1 (B_R \cap \Omega)} 
	\cap 
	{L_\infty (B_R \cap \Omega)}
$ 
is a solution of problem~\eqref{1.1} if
${A (x, D u)} \in {L_{p/(p-1)} (B_R \cap \Omega)}$,
$$
	- \int_\Omega
	A (x, D u)
	D \varphi
	\, dx
	+
	\int_\Omega
	b (x) |D u|^\alpha
	\varphi
	\, dx
	\ge
	0
$$
for any non-negative function
$
	\varphi 
	\in 
	C_0^\infty (B_R \cap \Omega),
$
and 
$
	u
	\psi 
	\in
	{
		\stackrel{\rm \scriptscriptstyle o}{W}\!\!{}_p^1
		(
			\Omega
		)
	}
$
for any 
$
	\psi
	\in
	C_0^\infty 
	(
		B_R
	).
$

In his classical papers~\cite{W1, W2}, N.~Wiener obtained a boundary point regularity criteria for solutions of the Dirichlet problem for the Laplace equation.
In other words, he found necessary and sufficient conditions for solutions of the Dirichlet problem for the Laplace equation to be continuous at a boundary point. The criteria was formulated in terms of capacity which is very similar to the one that arises in electrostatics. 
This approach proved to be very productive and was subsequently used by many authors~[1--7].
In paper~\cite{Mazya}, V.G.~Maz'ya managed to get sufficient regularity conditions for solutions of the Dirichlet problem for the p-Laplace equation. The results of V.G.~Maz'ya were generalized for quasilinear equations containing term with lower-order derivatives by R.~Gariepy and W.~Ziemer~\cite{GZ} and for systems of quasilinear equations by J.~Bj\"orn~\cite{Bjorn}. In so doing, authors of papers~\cite{Bjorn, GZ} imposed 
essential
restrictions on coefficients of the lower-order derivatives. 
In the case of problem~\eqref{1.1}, this restrictions take the form 
$b^{1 / (p - \alpha)} \in L_n (B_R \cap \Omega)$ if $1 < p < n$ 
and
$b^{1 / (p - \alpha)} \in L_\lambda (B_R \cap \Omega)$, $\lambda > n$, if $p = n$. 
Therefore, the results of~\cite{Bjorn, GZ} can not be applied if $b (x)$ grows fast enough as $x \to 0$ 
(see Examples~\ref{E2.1}--\ref{E2.3}).
Below we present Theorems~\ref{T2.1}--\ref{T2.10} that are free from this shortcoming. 

We use the following notations.
For every solution of~\eqref{1.1} we put
\begin{equation}
	M (r; u)
	=
	\operatorname*{ess\,sup}_{
		S_r
		\cap
		\Omega
	}
	u,
	\label{1.2}
\end{equation}
where the restriction of $u$ to
$
	S_r
	\cap
	\Omega,
$
$r \in (0, R)$,
is understood in the sense of the trace and the essential supremum in~\eqref{1.2} is taken with respect to $(n-1)$-dimensional Lebesgue measure on the sphere $S_r$.
In accordance with the maximum principle either ${M(\cdot; u)}$ is a monotonic function
on the whole interval $(0, R)$ or there exists $R_* \in (0, R)$ such that ${M(\cdot; u)}$
does not increase on $(0, R_*)$ and does not decrease on $(R_*, R)$.

Let $E$ be a non-empty open subset of the sphere $S_r$.
We denote
$$
	\lambda_{min} (E)
	=
	\inf_{
		\psi
		\in
		C_0^\infty (E)
	}
	\frac{
		\int_E
		|\nabla \psi|^p
		\,
		dS_r
	}{
		\int_E
		|\psi|^p
		\,
		dS_r
	},
$$
where
$
	|\nabla \psi|
	=
	(g^{ij} \nabla_i \psi \nabla_j \psi)^{1 / 2},
$
$g^{ij}$ is the dual metric tensor on $S_r$
induced by the standard euclidean metric on ${\mathbb R}^n$,
and $dS_r$ is the $(n-1)$-dimensional volume element of $S_r$.
By the variational principle, $\lambda_{min} (E)$ is the first eigenvalue of the problem
$$
	\Delta_p v
	=
	- \lambda
	|v|^{p-2} v
	\quad
	\mbox{in } E,
	\quad
	\left.
	v
	\right|_{\partial E}
	=
	0,
$$
for the $p$-Laplace-Beltrami operator
$
	\Delta_p v
	=
	\nabla_i
	(
		|\nabla v|^{p-2}
		g^{ij}
		\nabla_j v
	).
$

The capacity of a compact set $K \subset \omega$ relative to a non-empty open set 
$\omega \subset {\mathbb R}^n$ is defined as 
$$
    \operatorname{cap} (K, \omega)
    =
    \inf_\varphi \int_\omega
    |D \varphi|^p
    \,
    dx,
$$
where the infimum is taken over all functions 
$\varphi \in C_0^\infty (\omega)$ 
that are identically equal to one in a neighborhood of $K$.
By definition, the capacity of the empty set is equal to zero.
In the case of $\omega = {\mathbb R}^n$, we write $\operatorname{cap} (K)$ 
instead of $\operatorname{cap} (K, \omega)$.
If $p = 2$ and $n \ge 3$, then $\operatorname{cap} (K)$ coincides with the well-known Wiener capacity.

It can be shown that $\operatorname{cap} (K, \omega)$ has the following natural properties.

\begin{enumerate}
\item[(a)]{\it Monotonicity}:
If 
$K_1 \subset K_2$ 
and 
$\omega_2 \subset \omega_1$,
then
$$
	\operatorname{cap} (K_1,\omega_1) \le \operatorname{cap} (K_2,\omega_2).
$$
\item[(b)]{\it  Similarity property}:
If
$K' = \lambda K$ 
and 
$\omega' = \lambda \omega$, where $\lambda > 0$ is a real number, then
$$
	\operatorname{cap} (K',\omega') = \lambda^{n-p} \operatorname{cap} (K,\omega).
$$
\item[(c)]{\it Semiadditivity}:
Assume that $K_1$ and $K_2$ are compact subsets of an open set $\omega$, then
$$
	\operatorname{cap} (K_1 \cup K_2, \omega) 
	\le
	\operatorname{cap} (K_1,\omega)
	+
	\operatorname{cap} (K_2,\omega).
$$
\end{enumerate}

By the $\varepsilon$-essential inner diameter of an open set $\omega$, 
where $0 < \varepsilon < 1$ is a real number, we mean the value
$$
	\operatorname{diam}_\varepsilon \omega
	=
	\sup 
	\left\{
		r \in (0, \infty) 
		: 
		\exists x \in \omega 
		\;
		\frac{
			\operatorname{cap} (\overline{B_r^x} \setminus \omega, B_{2 r}^x)
		}{
			\operatorname{cap} (\overline{B_r}, B_{2 r})
		} 
		< 
		\varepsilon
	\right\}.
$$
In so doing, if $\omega = \emptyset$, then $\operatorname{diam}_\varepsilon \omega = 0$.

The $\varepsilon$-essential inner diameter is a monotone set function, i.e.
$
	\operatorname{diam}_\varepsilon \omega_1
	\le
	\operatorname{diam}_\varepsilon \omega_2
$
if $\omega_1 \subset \omega_2$.
It also is a monotone function of $\varepsilon$.
In other words, 
$
	\operatorname{diam}_{\varepsilon_1} \omega
	\le
	\operatorname{diam}_{\varepsilon_2} \omega
$
if $\varepsilon_1 \le \varepsilon_2$.

We say that $f \in {\mathcal L}_{\nu, \varepsilon} (\omega)$, 
where $\nu \ge 1$ and $0 < \varepsilon < 1$ are real numbers and $\omega$ is an open set, if 
$f \in {L_{\nu, loc} (\omega)}$ and 
$$
	\sup_{x \in \omega}
	\| f \|_{
		L_\nu (\omega \cap B_{\operatorname{diam}_\varepsilon \omega}^x)
	}
	<
	\infty.
$$
It can be seen that ${\mathcal L}_{\nu, \varepsilon} (\omega)$ is a Banach space with the norm 
$$
	\| f \|_{
		{\mathcal L}_{\nu, \varepsilon} (\omega)
	}
	=
	|S_1|^{- 1 / \nu}
	\sup_{x \in \omega}
	\| f \|_{
		L_\nu \left( \omega \cap B_{\operatorname{diam}_\varepsilon \omega}^x \right)
	},
$$
where $|S_1|$ is the $(n - 1)$-dimensional volume of $S_1$.
In the case of $f \in L_\infty (\omega)$, we obviously have
\begin{equation}	
\| f \|_{
		{\mathcal L}_{\nu, \varepsilon} (\omega)
	}
	\le
	(\operatorname{diam}_\varepsilon \omega)^{n / \nu}
	\| f \|_{
		L_\infty (\omega)
	}.
	\label{1.3}
\end{equation}

\section{Estimates of solutions near a boundary point}

Below we assume by default that $\Lambda$, $q$, and ${\mathcal D}$ are non-negative measurable functions such that
\begin{equation}
	\Lambda (r)
	\le
	\inf_{
		t
		\in 
		(r / \theta, r \theta)
		\cap
		(0, R)
	}
	\lambda_{min}(
		S_t
		\cap
		\Omega
	),
	\label{2.1}
\end{equation}
\begin{equation}
	q (r)
	\ge
	(\operatorname{diam}_\varepsilon \Omega_{r / \theta, r \theta})^{p - \alpha - n / \nu}
	\| b \|_{
		{\mathcal L}_{\nu, \varepsilon} (
			\Omega_{r / \theta, r \theta}
		)
	},
	\label{2.2}
\end{equation}
and
$$
	{\mathcal D} (r)
	\le
	\frac{
		1
	}{
		\operatorname{diam}_\varepsilon \Omega_{r / \theta, r \theta}
	}
	\label{2.6}
$$
for almost all $r \in (0, R)$, where
$\theta > 1$ and $0 < \varepsilon < 1$
are some real numbers.

\begin{Remark}\label{R2.1}
In view of~\eqref{1.3}, if $b \in L_\infty (\Omega_{r, R})$ for any $r \in (0, R)$,
then to perform~\eqref{2.2} it is sufficient to require that 
$$
	q (r)
	\ge
	(\operatorname{diam}_\varepsilon \Omega_{r / \theta, r \theta})^{p - \alpha}
	\operatorname*{esssup}_{
		\Omega_{r / \theta, r \theta}
	}
	b
$$
for almost all $r \in (0, R)$.
\end{Remark}

\begin{Theorem}\label{T2.1}
Let $p - 1 < \alpha \le p$ and
$$
	\int_0^R
	\frac{
		\min
		\{
			(r \Lambda (r))^{1 / (p - 1)},
			\Lambda^{1 / p} (r)
		\}
	}{
		1 + q^{1 / (\alpha - p + 1)} (r) 
	}
	\,
	d r
	=
	\infty.
$$
Then every non-negative solution of~\eqref{1.1} satisfies the estimate
\begin{equation}
	M (r; u)
	\le
	M (R; u)
	\exp
	\left(
		-
		C
		\int_r^R
		\frac{
			\min
			\{
				(t \Lambda (t))^{1 / (p - 1)},
				\Lambda^{1 / p} (t)
			\}
		}{
			1 + q^{1 / (\alpha - p + 1)} (t) 
		}
		\,
		dt
	\right)
	\label{T2.1.2}
\end{equation}
for all sufficiently small $r > 0$,
where the constant $C > 0$ depends only on
$n$, $p$, $\alpha$, $\varepsilon$, $\theta$, $\nu$, and the ellipticity constants $C_1$ and~$C_2$.
\end{Theorem}

\begin{Example}\label{E2.1}
Assume that $p - 1 < \alpha \le p$,
$
	\{
		(x', x_n) \in {\mathbb R}^n
		:
		|x'| < k_1 x_n,
		\:
		0 < x_n < R
	\}
	\subset
	B_R
	\setminus
	\Omega,
$
and
\begin{equation}
	b (x)
	\le
	k_2
	|x|^l
	\label{E2.1.1}
\end{equation}
for almost all $x \in B_R \cap \Omega$,
where $k_1$ and $k_2$ are positive constants and $l \in {\mathbb R}$.

If $l \ge \alpha - p$, then Theorem~\ref{T2.1} implies that
$M (r; u) \to 0$ as $r \to +0$
for any non-negative solution of~\eqref{1.1}.
In addition, the estimate
$$
	M (r; u)
	\le
	M (R; u)
	r^k
$$
is valid for all sufficiently small $r > 0$, where the constant $k > 0$ does not depend on $u$.
Really, we can take the function $\Lambda$ such that
\begin{equation}
	\Lambda (r) \sim r^{- p}
	\quad
	\mbox{as } r \to +0
	\label{E2.1.2}
\end{equation}
or, in other words, 
$$
	\varkappa_1
	r^{- p}
	\le
	\Lambda (r)
	\le
	\varkappa_2
	r^{- p}
$$
with some constants $\varkappa_1 > 0$ and $\varkappa_2 > 0$ for all $r > 0$ 
from a neighborhood of zero.
In so doing, as the $q$, we can take a bounded function.

We note that, from paper~\cite{GZ}, the required regularity follows only for $l > \alpha - p$. 
In the case of the critical exponent $l = \alpha - p$, the results of~\cite{GZ} are inapplicable.

Now, let the inequality
$$
	b (x)
	\le
	k_2
	|x|^{\alpha - p}
	\left(
		\log \frac{1}{|x|}
	\right)^\sigma
$$
be fulfilled instead of~\eqref{E2.1.1}.
In other words, we examine the case of the critical exponent $l = \alpha - p$. 
If $\sigma \le \alpha - p + 1$, then in accordance with
Theorem~\ref{T2.1}, where $\Lambda$ satisfies~\eqref{E2.1.2} and
$$
	q (r) 
	\sim 
	\left(
		\log \frac{1}{r}
	\right)^\sigma
	\quad
	\mbox{as } r \to +0,
$$
we have
$M (r; u) \to 0$ as $r \to +0$
for any non-negative solution of~\eqref{1.1}.
In addition, it can be shown that
\begin{equation}
	M (r; u)
	\le
	M (R; u)
	e^{- C f (r)}
	\label{E2.1.3}
\end{equation}
for all sufficiently small $r > 0$, where 
$$
	f (r)
	=
	\left\{
		\begin{aligned}
			&
			\log \frac{1}{r},
			&
			&
			\sigma \le 0,
			\\
			&
			\left(
				\log \frac{1}{r}
			\right)^{
				(\alpha - p + 1 - \sigma) / (\alpha - p + 1)
			},
			&
			&
			0 < \sigma < \alpha - p + 1,
			\\
			&
			\log \log \frac{1}{r},
			&
			&
			\sigma = \alpha - p + 1,
		\end{aligned}
	\right.
$$
and $C > 0$ is a constant independent of $u$.
\end{Example}

\begin{Theorem}\label{T2.2}
Let $p - 1 < \alpha \le p$, 
$$
	\int_0^R
	\frac{
		\Lambda^{1 / p} (r)
	}{
		1 + q^{1 / (\alpha - p + 1)} (r) 
	}
	\,
	d r
	=
	\infty
$$
and, moreover,
\begin{equation}
	\liminf_{r \to +0}
	r^{p - n}
	\operatorname{cap} 
	\left(
		\overline{B_{r \theta^{- 2 / 3}, r \theta^{- 1 / 3}}}
		\setminus
		\Omega,
		B_{r / \theta, r}
	\right)
	>
	0.
	\label{T2.2.1}
\end{equation}
Then every non-negative solution of~\eqref{1.1} satisfies the estimate
\begin{equation}
	M (r; u)
	\le
	M (R; u)
	\exp
	\left(
		-
		C
		\int_r^R
		\frac{
			\Lambda^{1 / p} (t)
		}{
			1 + q^{1 / (\alpha - p + 1)} (t) 
		}
		\,
		dt
	\right)
	\label{T2.2.2}
\end{equation}
for all sufficiently small $r > 0$,
where the constant $C > 0$ depends only on
$n$, $p$, $\alpha$, $\varepsilon$, $\theta$, $\nu$, $C_1$, $C_2$,
and on the limit in the left-hand side of~\eqref{T2.2.1}.
\end{Theorem}

\begin{Remark}\label{R2.2}
Condition~\eqref{T2.2.1} is obviously fulfilled if we can touch zero by a cone that lies entirely outside the set $\Omega$. This condition is also fulfilled if
$$
	\lim_{r \to +0}
	\frac{
		\operatorname{diam}_\varepsilon \Omega_{r / \theta, r \theta}
	}{
		r
	}
	=
	0.
$$
\end{Remark}

\begin{Theorem}\label{T2.3}
Let $p - 1 < \alpha \le p$,
$$
	\int_0^R
	\frac{
		{\mathcal D} (r)
	}{
		1 + q^{1 / (\alpha - p + 1)} (r) 
	}
	\,
	d r
	=
	\infty
$$
and, moreover,~\eqref{T2.2.1} holds.
Then every non-negative solution of~\eqref{1.1} satisfies the estimate
$$
	M (r; u)
	\le
	M (R; u)
	\exp
	\left(
		-
		C
		\int_r^R
		\frac{
			{\mathcal D} (t)
		}{
			1 + q^{1 / (\alpha - p + 1)} (t) 
		}
		\,
		dt
	\right)
$$
for all sufficiently small $r > 0$,
where the constant $C > 0$ depends only on
$n$, $p$, $\alpha$, $\varepsilon$, $\theta$, $\nu$, $C_1$, $C_2$,
and on the limit in the left-hand side of~\eqref{T2.2.1}.
\end{Theorem}

\begin{Example}\label{E2.2}
Assume that $p - 1 < \alpha \le p$, 
$
	B_R
	\cap
	\Omega
	\subset 
	\{
		(x', x_n) \in {\mathbb R}^n
		:
		|x_n| < k_1 |x'|^s,
		\:
		|x'| < R
	\}
$,
and~\eqref{E2.1.1} is valid,
where $k_1 > 0$, $k_2 > 0$, and $s > 1$ are some constants.

In the case of $l \ge \alpha - p + 1 - s$, applying Theorem~\ref{T2.3} with
$$
	{\mathcal D} (r)
	\sim
	r^{-s}
	\quad
	\mbox{and}
	\quad
	q (r)
	\sim
	r^{s (p - \alpha) + l}
	\quad
	\mbox{as } r \to +0,
$$
we obtain that $M (r; u) \to 0$ as $r \to +0$ for any non-negative solution of~\eqref{1.1}.
In so doing, Theorem~\ref{T2.3} implies estimate~\eqref{E2.1.3}, where
$$
	f (r)
	=
	\left\{
		\begin{aligned}
			&
			r^{1 - s},
			&
			&
			s (\alpha - p) \le l,
			\\
			&
			r^{
				(\alpha - p + 1 - s - l) / (\alpha - p + 1)
			},
			&
			&
			\alpha - p + 1 - s < l < s (\alpha - p),
			\\
			&
			\log \frac{1}{r},
			&
			&
			l = \alpha - p + 1 - s.
		\end{aligned}
	\right.
$$

We note that the results of paper~\cite{GZ} yields the required regularity for 
$l > (\alpha - p) (n + s - 1) / n$.
It does not present any particular problem to verify that 
$(\alpha - p) (n + s - 1) / n > \alpha - p + 1 - s$
for all positive integers $n$.
Thus, Theorem~\ref{T2.3} provides us with a regularity condition that is better than the analogous condition given in~\cite{GZ}.

Now, let the inequality
$$
	b (x)
	\le
	k_2
	|x|^{\alpha - p + 1 - s}
	\left(
		\log \frac{1}{|x|}
	\right)^\sigma
$$
be fulfilled instead of~\eqref{E2.1.1}.

If $\sigma \le \alpha - p + 1$, then
$M (r; u) \to 0$ as $r \to +0$
for any non-negative solution of~\eqref{1.1}.
In addition, the function ${M (\cdot; u)}$ satisfies estimate~\eqref{E2.1.3}, where 
$$
	f (r)
	=
	\left\{
		\begin{aligned}
			&
			\left(
				\log \frac{1}{r}
			\right)^{
				(\alpha - p + 1 - \sigma) / (\alpha - p + 1)
			},
			&
			&
			\sigma < \alpha - p + 1,
			\\
			&
			\log \log \frac{1}{r},
			&
			&
			\sigma = \alpha - p + 1.
		\end{aligned}
	\right.
$$
To show this, it is sufficient to apply Theorem~\ref{T2.3} with
$$
	{\mathcal D} (r)
	\sim
	r^{-s}
	\quad
	\mbox{and}
	\quad
	q (r)
	\sim
	r^{s (p - \alpha) + \alpha - p + 1 - s}
	\left(
		\log \frac{1}{r}
	\right)^\sigma
	\quad
	\mbox{as } r \to +0.
$$
\end{Example}

\begin{Theorem}\label{T2.4} 
Estimate~\eqref{T2.1.2} remains valid if, under the assumptions of Theorem~$\ref{T2.1}$,
the function $\Lambda$ satisfies the inequality
\begin{equation}
	\Lambda (r)
	\le
	\inf_{
		\Omega_{r \theta^{- 1 / 3}, r \theta^{1 / 3}}
	}
	\mu_\delta^p 
	+
	r^{-n}
	\operatorname{cap} 
	\left(
		\overline{B_{r \theta^{- 2 / 3}, r \theta^{- 1 / 3}}}
		\setminus
		\Omega,
		B_{r / \theta, r}
	\right)
	\label{T2.4.1}
\end{equation}
instead of~\eqref{2.1}, 
where $\theta > 1$ and $0 < \delta < 1 - \theta^{- 1 / 3}$ are some real numbers and 
$$
	\mu_\delta (x)
	=
	\sup_{
		r \in (0, \delta |x|)
	}
	(
		r^{1 - n}
		\operatorname{cap} (
			\overline{B_r^x}
			\setminus
			\Omega,
			B_{2 r}^x
		)
	)^{1 / (p - 1)}.
$$
In this case, the constant $C > 0$ in~\eqref{T2.1.2} depends also on $\delta$.
\end{Theorem}

\begin{Corollary}\label{C2.1}
Let the inequality
\begin{equation}
	\Lambda (r)
	\le
	r^{-n}
	\operatorname{cap} 
	\left(
		\overline{B_{r \theta^{- 2 / 3}, r \theta^{- 1/ 3}}}
		\setminus
		\Omega,
		B_{r / \theta, r}
	\right)
	\label{C2.1.1}
\end{equation}
be fulfilled instead of~\eqref{2.1} and, moreover,
$$
	\int_0^R
	\frac{
		(r \Lambda (r))^{1 / (p - 1)}
	}{
		1 + q^{1 / (\alpha - p + 1)} (r) 
	}
	\,
	d r
	=
	\infty.
$$
Then every non-negative solution of~\eqref{1.1} satisfies the estimate
$$
	M (r; u)
	\le
	M (R; u)
	\exp
	\left(
		-
		C
		\int_r^R
		\frac{
			(t \Lambda (t))^{1 / (p - 1)}
		}{
			1 + q^{1 / (\alpha - p + 1)} (t) 
		}
		\,
		dt
	\right)
$$
for all sufficiently small $r > 0$,
where the constant $C > 0$ depends only on
$n$, $p$, $\alpha$, $\theta$, $\nu$, and the ellipticity constants $C_1$ and~$C_2$.
\end{Corollary}

\begin{Theorem}\label{T2.5}
Estimate~\eqref{T2.2.2} remains valid if, in the assumptions of Theorem~$\ref{T2.2}$,
the function $\Lambda$ satisfies inequality~\eqref{T2.4.1} instead of~\eqref{2.1}.
In this case, the constant $C > 0$ in~\eqref{T2.2.2} depends also on $\delta$.
\end{Theorem}

\begin{Theorem}\label{T2.6}
Let $u$ be a non-negative solution of~\eqref{1.1}, where $\alpha = p - 1$.
Then there exist constants $k > 0$ and $C > 0$ depending only on
$n$, $p$, $\varepsilon$, $\theta$, $\nu$, and the ellipticity constants $C_1$ and $C_2$
such that the condition
\begin{equation}
	\int_0^R
	e^{- k q (r)}
	\min
	\{
		(r \Lambda (r))^{1/ (p - 1)},
		\Lambda^{1 /p} (r)
	\}
	\, 
	dr
	=
	\infty
	\label{T2.6.1}
\end{equation}
implies the estimate
\begin{equation}
	M (r; u)
	\le
	M (R; u)
	\exp
	\left(
		- C
		\int_r^R
		e^{- k q (t)}
		\min
		\{
			(r \Lambda (r))^{1/ (p - 1)},
			\Lambda^{1 /p} (r)
		\}
		\, 
		dt
	\right)
	\label{T2.6.2}
\end{equation}
for all sufficiently small $r > 0$.
\end{Theorem}

\begin{Corollary}\label{C2.2}
Let $u$ be a non-negative solution of~\eqref{1.1} with $\alpha = p - 1$
and, moreover,~\eqref{C2.1.1} holds instead of~\eqref{2.1}.
Then there exist constants $k > 0$ and $C > 0$ depending only on
$n$, $p$, $\alpha$, $\varepsilon$, $\theta$, $\nu$, and the ellipticity constants $C_1$ and $C_2$
such that the condition
$$
	\int_0^R
	e^{- k q (r)}
	(r \Lambda (r))^{1 / (p - 1)}
	\, 
	dr
	=
	\infty
$$
implies the estimate
$$
	M (r; u)
	\le
	M (R; u)
	\exp
	\left(
		- C
		\int_r^R
		e^{- k q (t)}
		(t \Lambda (t))^{1 / (p - 1)}
		\, 
		dt
	\right)
$$
for all sufficiently small $r > 0$.
\end{Corollary}

\begin{Theorem}\label{T2.7}
Let $u$ be a non-negative solution of~\eqref{1.1} with $\alpha = p - 1$
and, moreover,~\eqref{T2.2.1} holds.
Then there exist constants $k > 0$ and $C > 0$ depending only on
$n$, $p$, $\alpha$, $\varepsilon$, $\theta$, $\nu$, $C_1$, $C_2$,
and on the limit in the left-hand side of~\eqref{T2.2.1}
such that the condition
\begin{equation}
	\int_0^R
	e^{- k q (r)}
	\Lambda^{1 / p} (r)
	\, 
	dr
	=
	\infty
	\label{T2.7.1}
\end{equation}
implies the estimate
\begin{equation}
	M (r; u)
	\le
	M (R; u)
	\exp
	\left(
		- C
		\int_r^R
		e^{- k q (t)}
		\Lambda^{1 / p} (t)
		\, 
		dt
	\right)
	\label{T2.7.2}
\end{equation}
for all sufficiently small $r > 0$.
\end{Theorem}

\begin{Theorem}\label{T2.8}
Let $u$ be a non-negative solution of~\eqref{1.1} with $\alpha = p - 1$
and, moreover,~\eqref{T2.2.1} holds.
Then there exist constants $k > 0$ and $C > 0$ depending only on
$n$, $p$, $\alpha$, $\varepsilon$, $\theta$, $\nu$, $C_1$, $C_2$,
and on the limit in the left-hand side of~\eqref{T2.2.1}
such that the condition
$$
	\int_0^R
	e^{- k q (r)}
	{\mathcal D} (r)
	\, 
	dr
	=
	\infty
$$
implies the estimate
$$
	M (r; u)
	\le
	M (R; u)
	\exp
	\left(
		- C
		\int_r^R
		e^{- k q (t)}
		{\mathcal D} (t)
		\, 
		dt
	\right)
$$
for all sufficiently small $r > 0$.
\end{Theorem}

\begin{Example}\label{E2.3}
Assume that $\alpha = p - 1$, 
$
	B_R
	\cap
	\Omega
	\subset 
	\{
		(x', x_n) \in {\mathbb R}^n
		:
		|x_n| < k_1 |x'|^s,
		\:
		|x'| < R
	\}
$,
and~\eqref{E2.1.1} is valid,
where $k_1 > 0$, $k_2 > 0$, and $s > 1$ are some constants.

In the case of $l \ge - s$, taking
$$
	{\mathcal D} (r)
	\sim
	r^{-s}
	\quad
	\mbox{and}
	\quad
	q (r)
	\sim
	1
	\quad
	\mbox{as } r \to +0
$$
in Theorem~\ref{T2.8}, we obtain $M (r; u) \to 0$ as $r \to +0$ for any non-negative solution of~\eqref{1.1}. In so doing, estimate~\eqref{E2.1.3} is valid, where
$$
	f (r)
	=
	r^{1 - s}.
$$

Note that the results of paper~\cite{GZ} guarantee the required regularity for 
$l > - (n + s - 1) / n$.
It is easy to see that $- s < - (n + s - 1) / n$ for all integers $n \ge 2$.
Thus, Theorem~\ref{T2.8} gives us a better regularity condition than the results of paper~\cite{GZ}.
\end{Example}

\begin{Theorem}\label{T2.9}
In the hypotheses of Theorem~$\ref{T2.6}$, let the function $\Lambda$ satisfies inequality~\eqref{T2.4.1} instead of~\eqref{2.1}. 
Then there exist constants $k > 0$ and $C > 0$ depending only on
$n$, $p$, $\delta$, $\varepsilon$, $\theta$, $\nu$, and the ellipticity constants $C_1$ and~$C_2$
such that the condition~\eqref{T2.6.1} implies estimate~\eqref{T2.6.2}.
\end{Theorem}

\begin{Theorem}\label{T2.10}
In the hypotheses of Theorem~$\ref{T2.7}$, let the function $\Lambda$ satisfies inequality~\eqref{T2.4.1} instead of~\eqref{2.1}. 
Then there exist constants $k > 0$ and $C > 0$ depending only on
$n$, $p$, $\delta$, $\alpha$, $\varepsilon$, $\theta$, $\nu$, $C_1$, $C_2$,
and on the limit in the left-hand side of~\eqref{T2.2.1}
such that the condition~\eqref{T2.7.1} implies estimate~\eqref{T2.7.2}.
\end{Theorem}


\end{document}